\input amstex
\magnification=1200
\documentstyle{amsppt}
\NoBlackBoxes
\NoRunningHeads
\topmatter
\title Interactive forces\endtitle
\author Denis V. Juriev\endauthor
\affil ul.Miklukho-Maklaya 20-180, Moscow 117437 Russia\linebreak
(e-mail: denis\@juriev.msk.ru)\endaffil
\abstract\nofrills This short article is devoted to the dynamics of controlled
(and, therefore, open) systems. The internal forces, which appear only in the 
presence of external free controls and depend explicitely on them, are 
considered. Such {\it interactive forces} may be regarded as feedbacks 
generated in the system by the external free controls. In particular, one is 
able to interpret the interactive controls as couplings of free controls with 
the action of interactive forces; such dynamical approach to the interactive 
systems often simplifies their analysis and allows to perform their synthesis 
on a systematic level. The field theoretic counterparts of interactive forces, 
the fields of interactive forces, are also considered.
\endabstract
\endtopmatter
\document
The original impetus to write this article was an intention to construct a
bridge between the note [1] on the physical dynamics of information processes 
and the series [2] on the mathematics of interactive games. However, it was 
understood soon that the subject is more general and should be classified as
dynamics of controlled (and, therefore, open) systems as a part of 
mathematical physics. Such approach enlarges the area of researches, makes
the understanding deeper and physically pithy as well as allows to adopt
simple mechanical models for analysis of interactive systems.

The new proposed concept of the dynamics of controlled system is one of
the interactive forces. Interactive forces are internal forces in the
controlled systems, which appear only in the presence of external free
controls and explicitely depends on them. In classical mechanics with
constraints one is able to treat the reactions of constraints as interactive
forces if the external controls are realized as forces themselves. This
example is simplest and other ones may be regarded as its sophisticated
generalizations. One can look for interactive forces in mathematical
psychology, mathematical economics, mathematical sociology and other 
behavioral sciences, where the nontrivial reactions of the systems on the 
external controls are considered. The dynamical approach to the interactive 
systems, the controlled systems, in which free controls are coupled with the 
unknown or incompletely known feedbacks, thus interpreted as an action of the 
interactive forces, symplifies their analysis and allows to perform their 
synthesis on a systematic level.

The field theoretic counterparts of interactive forces, the fields of 
interactive forces, are also considered.

Of course, the mechanical approach to the interactive forces does not
substitute the physical analysis of their nature, which may be rather
different even for mechanically similar systems. Also such approach has its
evident methodological boundaries when one has deal with interactive systems, 
which may be regarded as complex systems with incomplete information [2], and
is often practically useful only together with other approaches. However,
the explication of this new ``coordinate axis'' for the description of 
interactive systems should not be underestimated.

\subhead 1. Interactive forces\endsubhead

\definition{Definition 1} {\it Interactive forces\/} are internal forces in 
the controlled systems, which appear only in the presence of external free
controls and explicitely depends on them.
\enddefinition

\remark{Example 1} In classical mechanics with constraints one is able to 
treat the reactions of constraints as interactive forces if the external 
controls are realized as forces themselves.
\endremark

\remark{Example 2 (A quantum mechanical example)} In quantum mechanics of
continuously observed particle one may treat the Belavkin counterterm in
quantum Hamiltonian [3] as corresponding to the interactive forces.
\endremark

\remark{Remark 1} One is able to consider the fields of interactive forces
in the Faradey-Maxwell sense.
\endremark

\subhead 2. Interactive systems\endsubhead 
The concept of an interactive control was recently proposed by the author 
[2] to take into account the complex composition of controls of a real human 
person, which are often complicated couplings of his/her cognitive and known 
controls with the unknown subconscious behavioral reactions. This formalism 
is applicable also to the description of external unknown influences and, 
thus, is useful for problems in computer science (e.g. the semi-artificially 
controlled distribution of resources), mathematical economics (e.g. the 
financial games with unknown dynamical factors) and sociology (e.g. the 
collective decision making).

\definition{Definition 2 [2]} An {\it interactive system\/} (with $n$
{\it interactive controls\/}) is a control system with $n$ independent 
controls coupled with unknown or incompletely known feedbacks (the feedbacks
as well as their couplings with controls are of a so complicated nature that 
their can not be described completely).
\enddefinition

Below we shall consider only deterministic and differential interactive
systems. In this case the general interactive system may be written in the 
form:
$$\dot\varphi=\Phi(\varphi,u_1,u_2,\ldots,u_n),$$
where $\varphi$ characterizes the state of the system and $u_i$ are
the interactive controls:
$$u_i(t)=u_i(u_i^\circ(t),\left.[\varphi(\tau)]\right|_{\tau\leqslant t}),$$
i.e. the independent controls $u_i^\circ(t)$ coupled with the feedbacks on
$\left.[\varphi(\tau)]\right|_{\tau\leqslant t}$. One may suppose that the
feedbacks are integrodifferential on $t$.

However, it is reasonable to consider the {\it differential interactive
systems}, whose feedbacks are purely differential. It means that
$$u_i(t)=u_i(u_i^\circ(t),\varphi(t),\ldots,\varphi^{(k)}(t)).$$
A reduction of general interactive systems to the differential ones via the
introducing of the so-called {\it intention fields\/} was described in [2]. 
Such intention fields may be considered as hidden parameters. They are often
considered in the non-behaviourist psychology on the qualitative level. One 
is able to interpret them as internal (individual) as well as external (e.g. 
collective) parameters. 

The interactive systems introduced above may be generalized in the following
way, which leads to the {\it indeterminate interactive systems},
is based on the idea that the pure controls $u_i^\circ(t)$ and the 
interactive controls $u_i(t)$ should not be obligatory related in the
considered way. More generally one should only postulate that there are
some time-independent quantities $F_\alpha(u_i(t),u_i^\circ(t),\varphi(t),
\ldots,\varphi^{(k)}(t))$ for the independent magnitudes $u_i(t)$ and 
$u_i^\circ(t)$. Such claim is evidently weaker than one of Def.2. For
instance, one may consider the inverse dependence of the pure and 
interactive controls: $u_i^\circ(t)=u_i^\circ(u_i(t),\varphi(t),\ldots,
\varphi^{(k)}(t))$.

The inverse dependence of the pure and interactive controls has a nice
psychological interpretation. Instead of thinking of our action consisting
of the conscious and unconscious parts and interpreting the least as 
unknown feedbacks which ``dress'' the first, one is able to consider
our action as a single whole whereas the act of consciousness is in
the extraction of a part which it declares as its ``property''. So
interpreted the function $u_i^\circ(u_i,\varphi,\ldots,\varphi^{(k)})$ 
realizes the ``filtering'' procedure.

\subhead 3. Interactive forces in interactive systems\endsubhead
If the pure controls $u^\circ_i(t)$ in an interactive systems have a dynamical 
character, i.e. are realized by the external forces, one is able to try to 
interpret the action of interactive controls $u_i(t)$ as a sum of an action of 
pure controls and interactive forces $F_\alpha(u^\circ_i(t),\varphi(t))$. 
Precisely, let the system has the generalized Newtonian form
$$\left\{\aligned
\dot\bold q&=\tfrac1m\bold p\\
\dot\bold p&=F(\vec u,\bold q,\bold p)
\endaligned\right.$$
where $\bold q$ and $\bold p$ are generalized coordinates and momenta,
respectively, and $\vec u=(u_1,\ldots u_n)$ is the vector of $n$ independent
interactive controls. Then this system may be represented as
$$\left\{\aligned
\dot\bold q&=\tfrac1m\bold p\\
\dot\bold p&=F(\vec u^\circ,\bold q,\bold p)+\sum_\alpha F_\alpha(\vec u^\circ,
\bold q,\bold p)
\endaligned\right.$$
where $\vec u^\circ=(u_1^\circ,\ldots,u_n^\circ)$ is the vector of $n$ 
independent pure controls and $F_\alpha$ are the interactive forces.

If the forces are potential one may consider dynamics in Hamiltonian form
instead of the generalized Newtonian one. In this case the Hamiltonian 
$H(\vec u,\bold p,\bold q)$ is decomposed into the sum $H(\vec u^\circ,\bold 
p,\bold q)+\sum_\alpha H_\alpha(u^\circ,\bold p,\bold q)$, where $H_\alpha$ 
are the terms corresponding to the interactive forces.

This dynamical interpretation of interactive controls often symplifies an
analysis of the interactive systems and allows to perform their synthesis
on a systematic level. For instance, one is able to apply various mechanical
principles such as superposition principle or the third law of dynamics to
the interactive systems. Moreover, and it should be specially emphasized that 
{\sl the synthesis of interactive systems is straightforward if their 
dynamical nature is explicated}. This is just the same as for the classical 
mechanics, where the Hamiltonians for two isolated systems do not provide a 
sufficient information to reconstruct the Hamiltonian of interaction, whereas 
the knowledge of dynamical nature of both systems completely determines their 
interaction.

\remark{Remark 2 (The actions of interactive forces)} An action of the pure 
control may be zero in the mechanical sense (i.e. it does not change the
energy of system) or, otherwords, the pure control ``does not work''
(for instance, the corresponding force is magnetic-type of gyroscopic), 
whereas the interactive forces are acting and, therefore, change the energy 
of system. In this situation the whole action is simulateously effective 
and energy-expenseless for the controller so one may say that the pure control 
is resulted only in such non-energetical {\sl effort\/} whereas the interactive 
forces really {\sl act}.
\endremark

\remark{Remark 3} Note that the concept of a field of interactive forces is
close to the concept of an intention field. However, one is able to consider
the intention fields, which are not defined by fields of interactive forces.
And there are fields of interactive forces (see examples 1,2 above), which
are not constructed as intention fields in any interactive systems.
\endremark

\remark{Remark 4} One may treat intention fields and fields of interactive
forces in physical interactive information systems as certain 
``objectivizations'' of {\sl m\B ay\B a}, the philosophical concept adapted 
in [1] for a description of dynamics of physical interactive information 
systems. So the intention fields and fields of interactive forces are 
mathematical and dynamical counterparts for each other in this particular 
case.
\endremark

It is rather interesting to explicate the dynamical aspects of the perception
processes and image understanding, which interactive game theoretic 
description was done in [4]. Indeed, any interior image in the subjective
perceptional space of observer is a configuration of fields of interactive
forces due to its self-generating character. Moreover, any exterior object
also should be described in terms of such fields in view of its dynamical
property to generate the interior images. This is certainly reminiscent of 
the old Ernst Mach's ideas.

As it was mentioned in the preceeding articles the visualization of intention
fields is an effective tool for the organization of interactive processes.
The visualization of such fields, which correspond to the fields of 
interactive forces, explicates the dynamical aspect of this procedure.
This allows to use the dynamical nature of visual perception as a tool for 
the visualization of interactive processes of another nature. On the other
hand, an application of visualization to the visual perception itself
may be regarded as a {\sl dynamical exteriorization\/} of internal visual
perception processes and, thus, is effective for their description and the 
understanding of their dynamical nature.

\Refs
\roster
\item"[1]" Juriev D., On the dynamics of physical interactive information 
systems. E-print: mp\_arc/97-158.
\item"[2]" Juriev D., Interactive games and representation theory. I, II.
E-prints: math.FA/980320, math.RT/9808098; Games, predictions, interactivity. 
E-print: math.OC/9906107.
\item"[3]" Belavkin V.P., Phys. Lett. 144A (1989) 355; Belavkin V., Staszewski
P., Phys. Rev. 45A (1992) 1347; Belavkin V.P., Kolokoltsov V.N., Teor.Matem.
Fiz. 89 (1991) 163; Kolokoltsov V.N., Matem.Zametki 50 (1991).
\item"[4]" Juriev D., Perception games, the image understanding and 
interpretational geometry. E-print: math.OC/9905165; Kaleidoscope-roulette:
the resonance phenomena in perception games. E-print: math.OC/9905180.
\endroster
\endRefs
\enddocument